\documentclass[10pt]{article}

\usepackage{amsmath,amsthm,amssymb,xurl}
\usepackage[hidelinks]{hyperref}
\usepackage{tikz,authblk}
\usetikzlibrary{shapes.geometric}

\renewcommand{\Pr}{{\ensuremath{\mathrm{Pr}}}}
\allowdisplaybreaks

\title{Daniel Litt's Probability Puzzle}
\author[1]{Maura B.\ Paterson}
\affil[1]{School of Computing and Mathematical Sciences, Birkbeck, University of London, Malet St, London WC1E 7HX, UK}
\author[2]{Douglas R.\ Stinson\thanks{D.R.\ Stinson's research is supported by  NSERC discovery grant RGPIN-03882.}}
\affil[2]{David R.\ Cheriton School of Computer Science\\University of Waterloo\\ Waterloo ON, N2L 3G1, Canada}

\begin{document}
\maketitle

\begin{abstract}
In this expository note, we discuss a ``balls-and-urns'' probability puzzle posed by Daniel Litt.
\end{abstract}

\section{Introduction}

On August 29, 2024, Quanta Magazine published a discussion \cite{quanta}
of a probability puzzle that Daniel Litt posed on X on January 29, 2024. Here is the description of the puzzle from \cite{quanta}:
\begin{quote}
Imagine, he wrote, that you have an urn filled with 100 balls, some red and some green. You can’t see inside; all you know is that someone determined the number of red balls by picking a number between zero and 100 from a hat. You reach into the urn and pull out a ball. It’s red. If you now pull out a second ball, is it more likely to be red or green (or are the two colours equally likely)?
\end{quote}
These kinds of probability puzzles can be notoriously tricky, and intuition can easily lead one astray. For these reasons, we feel it is often useful to approach such puzzles from a mathematical point of view. We present three rigorous approaches to solving the problem in the next section. In Section \ref{inductive.sec}, we use an inductive approach to solve the puzzle. 
Finally, Section
\ref{comments.sec} concludes with some comments and discussion.

\section{Three solutions}

\subsection{Probability distribution on the urns}
\label{sec1.sec}

For $1 \leq x \leq 100$, let $U_x$ denote an urn that contains $x$ red balls and $100-x$ green balls. Suppose  a red ball has been removed from $U_x$.  If we remove a second ball from $U_x$, the probability this second ball is red is
$(x-1)/99$. Let's denote this (conditional) probability by $P_{x}$. The overall probability $P$ that the second ball is red depends on the probability distribution on the set of possible urns, $\mathcal{U} = \{U_1, \dots , U_{100}\}$.

Suppose we (incorrectly) assume a uniform distribution on $\mathcal{U}$, i.e., $\Pr [ U_x ] = 1/100$ for all $x$, $1 \leq x \leq 100$. In this case, 
\begin{align}
P &= \sum _{x = 1}^{100} \Pr [ U_x ] \times P_x\nonumber\\
&= \sum _{x = 1}^{100} \left( \frac{1}{100} \times \frac{x-1}{99} \right)\label{sum0.eq}\\
&= \frac{1}{100(99)} \times \sum _{x = 1}^{100} (x-1)\nonumber\\
&= \frac{1}{100(99)} \times \frac{99(100)}{2}\nonumber\\
&= \frac{1}{2}.\nonumber
\end{align}
This answer would be correct in the scenario where $U_x$ is first chosen uniformly  at random and then a red ball is removed from $U_x$. However, this is not quite the scenario that is described in the puzzle. 

The puzzle  assumes that urns with more red balls are more likely to be chosen than urns with fewer red balls. Imagine that we have created 100 urns $U_1, \dots , U_{100}$, as previously described.
The total number of red balls in all the urns is 
\[ \sum _{x = 1}^{100} x = \frac{100(101)}{2}.\]
If we choose a red ball uniformly at random, then the probability that the ball is chosen from urn $U_x$ is \[
\Pr [ U_x ] = \frac{x}{\frac{100(101)}{2}} = \frac{2x}{100(101)}.\]
Now we compute 
\begin{align}
P &= \sum _{x = 1}^{100} \Pr [ U_x ] \times P_x\nonumber\\
&= \sum _{x = 1}^{100} \frac{2x}{100(101)} \times \frac{x-1}{99}\nonumber\\
&= \frac{2}{99(100)(101)} \times \sum _{x = 1}^{100} (x^2-x)\label{sum.eq}\\
&= \frac{2}{99(100)(101)} \times \left( \frac{100(101)(201)}{6} - \frac{100(101)}{2} \right)\nonumber\\
&= \frac{2}{99} \left( \frac{201}{6} - \frac{1}{2} \right)\nonumber\\
&= \frac{2}{99} \left( \frac{198}{6}\right)\nonumber\\
&= \frac{2}{3}.\nonumber
\end{align}
The value $2/3$ is in fact the intended answer to the question. The reason the puzzle is rather tricky is possibly due to the nonuniform distribution on $\mathcal{U}$ that is induced by the initial choice of a red ball.

\bigskip

\noindent{\bf Remarks:}
\begin{enumerate}
\item When we compute $\sum x^2$ in (\ref{sum.eq}), we obtain a \emph{square pyramidal number}; see \cite{SPN}. The study of these numbers dates back to Archimedes. The classical formula  $\sum_{x=1}^n x^2 = n(n+1)(2n+1)/6$ can be interpreted as the volume of a square pyramid where the levels consist of square grids of $x^2$ unit spheres, $x = 1, \dots , n$.
\item
We also observe that computation in (\ref{sum.eq}) can be simplified slightly by using the \emph{hockey stick identity}
\[ \sum_{x = r}^{n} \binom{x}{r} = \binom{n+1}{r+1},\]
 setting $r = 2$ and $n = 100$.

\end{enumerate}

\subsection{A solution based on conditional probability}

Consider the following formulation of the puzzle. Suppose that $x \in \{0, \dots , 100\}$ is chosen uniformly at random. Then $x$ red balls and $100-x$ green balls are placed in an urn. Finally, two balls are drawn without replacement from the urn. We are asked to compute the conditional probability \[ P = \Pr[ \text{both balls are red} \mid \text{the first ball is red}] =
\frac{\Pr[ \text{both balls are red}]}{\Pr[\text{the first ball is red}]} .\]
It is clear that \[ \Pr[ \text{the first ball is red} ] = \frac{1}{2},\] where the probability is computed over the set of all $101$ possible urns.
To compute $\Pr[ \text{both balls are red} ]$, we observe that
\[ \Pr[ \text{the urn contains $x$ red balls and both chosen balls are red} ]
= \frac{1}{101} \times \frac{x(x-1)}{100(99)}.\]
Hence,  
\[ \Pr[ \text{both balls are red}] = \sum_{x=0}^{100} \frac{x(x-1)}{101(100)(99)}.\]
Therefore, 
\[ P = \frac{2}{101(100)(99)}\sum_{x=0}^{100} x(x-1).\]
Thus $P = 2/3$, by the same calculation that was performed in Section \ref{sec1.sec}.

\subsection{A symmetry argument}

The article \cite{quanta} also presents the following elegant argument (attributed to George Lowther in a post on ``X'') that the answer is $2/3$. This argument involves almost no mathematical calculations: 
\begin{quote}
Imagine that instead of starting with 100 balls, you start with 101 balls in a row. Pick a ball at random. Then colour the balls to the left of it green and the ones to the right of it red. Throw that ball away, leaving 100 balls.

Then pick a second ball at random. That ball corresponds to the first ball in the original problem. The problem tells you that you picked a red ball, so it was to the right of the ball you threw away. Now pick a third ball. This ball is either to the left of the first ball, between the first ball and the second, or to the right of the second. In two of the three possibilities, the third ball is red. So the probability that the ball is red is 2/3.
\end{quote}
George Lowther says in his ``X'' post that the answer is $2/3$ ``by symmetry.'' This certainly seems reasonable, but it might be of interest to develop a more detailed mathematical symmetry-based argument. In other words, how do we formalize this notion of symmetry?

We consider a slightly different, but equivalent, formulation of this argument. Suppose we have the following three-step process, starting from 100 balls (instead of starting from 101 balls and throwing the first ball away).
\begin{enumerate}
\item Choose an index $i$ uniformly at random from the set $\{0, \dots , 100\}$. The first $i$ balls are coloured green and the remaining $100-i$ balls are coloured red.
\item Choose an index $j$ uniformly at random from the set $\{i+1, \dots , 100\}$. This is the first ball chosen, which is red. 
\item Choose an index $k$ uniformly at random from the set $\{1, \dots , 100\} \setminus \{j\}$. This is the second ball chosen, which can be either red or green.
\end{enumerate}
Let's analyse steps 1 and 2. There are $\binom{100}{2}$ possible choices for the ordered pair $(i,j)$, because
$0 \leq i < j \leq 100$. For each choice of $i$ and $j$, we obtain three intervals (some of which may be empty): $[1,i]$, $[i+1,j-1]$ and $[j+1,100]$. The first  interval corresponds to green balls and the last two intervals correspond to red balls. Note that we are ignoring the red ball chosen in step 2, which occurs in the $j$th position.

In step 3, the chosen ball is red if and only if $k$ is in the second or third interval.
The sum of the lengths of all three intervals equals $99$. So we can express the lengths of the three intervals as a triple $(\ell_1,\ell_2,\ell_3)$ of non-negative integers whose sum is 99. As $i$ and $j$ vary over all $\binom{100}{2}$ possibilities, we get all possible such triples occurring once each. It then follows immediately (by symmetry!) that the expectation of each $\ell_i$ is 33. By linearity of expectation, the expected value of $\ell_2 + \ell_3$ is 66. (Alternatively, we can observe that $\ell_2 + \ell_3 = 99 - \ell_1$. Since the expectation of $\ell_1$ is $33$, the expectation of $99 - \ell_1$ must be $66$.) So the probability that the second ball is red is $66/99 = 2/3$. 

\section{An inductive approach}
\label{inductive.sec}

We could consider the generalisation of the puzzle where we start with an arbitrary fixed number $n \geq 2$ of balls.
The puzzle always has the same solution, independent of the value of $n$, namely $2/3$. It suffices to looks at any of the proofs we have presented and replace ``100'' by ``$n$''; the proofs are otherwise unchanged.
However, if we are thinking in terms of an arbitrary number of balls, we could consider an inductive approach to analysing the problem. It turns out that an inductive approach uses less computation---similar to the symmetry argument---in the sense that we do not have to evaluate sums of squares of integers.

Since we will be considering varying sizes of urns, let us denote by $U_{n,x}$ an urn  that contains $n$ balls, $x$ of which are red (where $0 \leq x \leq n$). Fix a value of $n$ and choose a random value $x \in \{0, \dots ,n\}$. Suppose we sample two balls from $U_{n,x}$ without replacement. A sample will be denoted as one of three possible multisets $\{R,R\}$, 
$\{R,G\}$ or $\{G,G\}$. For $j = 0, 1,2$, let $\mathbf{E}_j$ denote the event that the sample contains $j$ red balls (and hence $2-j$ green balls). We will prove by induction on $n$ that 
$\Pr [\mathbf{E}_j] = 1/3$, for $j = 0,1,2$.

The base case (where $n = 2$) is clearly true, since there is only one possible sample from each $U_{2,x}$ (for $x = 0,1,2$) and the three samples are respectively $\{G,G\}$, $\{R,G\}$ and $\{R,R\}$.

Now, as an inductive hypothesis, assume that $\Pr [\mathbf{E}_j] = 1/3$ ($j = 0,1,2$), for urns containing $n$ balls. 
For $0 \leq x \leq n$, we can construct $U_{n+1,x}$ from $U_{n,x}$ by adding a new green ball to the urn. The urn 
$U_{n+1,n+1}$ contains $n+1$ new red balls.

Let us consider the possible samples of two balls from a randomly chosen $U_{n+1,x}$. If we restrict to the samples consisting of two ``old'' balls, we have 
$\Pr [\mathbf{E}_j] = 1/3$, for $j = 0,1,2$, by induction. Let us now examine samples that contain at least one new ball. We observe the following:
\begin{itemize}
\item There are $n(n+1)/2$ samples that contain a new green ball and an old green ball 
and there are $n(n+1)/2$ samples that contain a new green ball and an old red ball. This follows because the urns
$U_0, \dots , U_n$ contain an equal number (namely, $n(n+1)/2$) of red and green balls.
\item There are $n(n+1)/2$ samples that contain two new red balls (these are the samples from $U_{n+1,n+1}$).
\end{itemize}
Therefore, among the samples that contain at least one new ball, we have $\Pr [\mathbf{E}_j] = 1/3$, for $j = 0,1,2$.
Then it follows (by induction) that $\Pr [\mathbf{E}_j] = 1/3$, for all $n \geq 2$ and for  $j = 0,1,2$.

With this result in hand, it is a simple matter to solve Litt's puzzle. In Litt's puzzle, we choose two balls in succession. Each unordered sample of two balls ($\{G,G\}$, 
$\{R,G\}$ or $\{R,R\}$) can be ordered in two equally probably ways. The corresponding  outcomes are 
$(G,G)$, $(G,G)$, $(G,R)$, $(R,G)$, $(R,R)$ and $(R,R)$.
Therefore we have
\[
\Pr[(G,G)] = \frac{1}{3}, \quad
\Pr[(R,R)] = \frac{1}{3}, \quad
\Pr[(R,G)] = \frac{1}{6}, \quad
\Pr[(G,R)] = \frac{1}{6}.
\]
We want to compute
\[ P = \Pr[ \text{both balls are red} \mid \text{the first ball is red}].\]
We see immediately that
\[ P = \frac{\frac{1}{3}}{\frac{1}{3} + \frac{1}{6}} = \frac{\frac{1}{3}}{\frac{1}{2}} = \frac{2}{3}.\]

\noindent{\bf Remark:}
The computation done to evaluate (\ref{sum0.eq}) has a nice interpretation in this inductive setting.
If we remove the urn $U_{n,0}$ and then remove one red ball from each of remaining urns, we end up with the urns
$U_{n-1,x}$ for $0 \leq x \leq n-1$ (in the case of  (\ref{sum0.eq}), we have $n = 100)$. 
The calculation in (\ref{sum0.eq}) just computes the probability that a randomly chosen ball in the collection of urns $U_{99,x}$ is red; this probability is of course just $1/2$.

\section{Comments}
\label{comments.sec}

The $n=2$ case of Litt's puzzle is a classical puzzle, known as ``Bertrand's box paradox.'' It dates back to 1889, where it was described in Joseph Bertrand's book \emph{Calcul des Probabilit\'{e}s}; see \cite{BBP}.
A perhaps better-known puzzle of this general type is the ``boy or girl paradox'' that is discussed in \cite{boyorgirl}. This puzzle was posed by Martin Gardner in 1959 in his long-running Mathematical Games column in \emph{Scientific American}.
The ``paradox'' is a consequence of the following two questions, as stated in \cite{boyorgirl}:
\begin{quote}
\begin{enumerate} \item Mr. Jones has two children. The older child is a girl. What is the probability that both children are girls?
\item Mr. Smith has two children. At least one of them is a boy. What is the probability that both children are boys?
\end{enumerate}
\end{quote}
The answer to the first question is clearly $1/2$, if we assume that the gender of the younger child is independent of the gender of the first child.

The second question is a bit trickier, and it has been argued that the answer is $1/2$ as well as $1/3$. The answer of $1/3$ is obtained as follows. We consider the four possible distributions of two children to be equally likely:
$BB$, $BG$, $GB$, $GG$. If at least one child is a boy, then $GG$ is impossible, so there remain three possible distributions. In two of these three distributions---namely $BG$ and $GB$---the ``other'' child is a girl. Therefore the answer is $1/3$. 

\medskip

Of course, Daniel Litt's puzzle has the answer $2/3$, which is different from either version of the boy or girl paradox. Basically, everything comes down to the underlying distribution. In the two-ball version of Litt's puzzle (i.e., Bertrand's box paradox), we have discussed how the underlying distributions of balls in urns yield three equally likely (unordered) multisets: $\{ R,R\}$, $\{R,G\}$, and $\{G,G\}$. So there is a uniform distribution on the possible multisets of cardinality two. (This was the base case of the inductive approach we discussed in Section \ref{inductive.sec}. We also showed that the same distribution holds, independent of the number of balls.) In contrast, in the boy or girl paradox, there is a uniform distribution defined on ordered pairs (of children).

\medskip

In many probability puzzles, the challenge is to convert an English language description of a problem into mathematical language. It is also important to be cognizant of any implicit assumptions that are being made. For example, in Litt's puzzle, it is implicitly assumed that all the balls in an urn are equally likely to be chosen. However, if we view the urn as a physical object and we think of reaching into the urn to choose a ball, then this assumption may not be valid.  Someone who places the balls into the urn might place all the green balls near the top, so they are more likely to be chosen.

\section*{Acknowledgements}

We thank Donald Kreher and Bill Martin for helpful discussions.

\end{document}